\documentclass{amsart}
\usepackage{amsfonts}
\input{epsf}
\newtheorem{theorem}{Theorem}[section]
\newtheorem{lemma}[theorem]{Lemma}

\theoremstyle{definition}

\theoremstyle{remark}

\numberwithin{equation}{section}

\newcommand{\ee}{{\bf E}}
\newcommand{\R}{\mathbb{R}}
\newcommand{\myprob}{{{\mathcal M}_1}}

\begin{document}

\title{A Reversion of the Chernoff Bound}

\author{Ted Theodosopoulos}
\address{Department of Decision Sciences and \\ Department of Mathematics 
\\ Drexel University \\ Philadelphia, PA 19066}
\email{theo@drexel.edu}
\urladdr{www.lebow.drexel.edu/theodosopoulos}


\date{January 21, 2005.}



\begin{abstract}
This paper describes the construction of a lower bound for the tails of 
general random variables, using solely knowledge of their moment generating 
function.  The {\it tilting} procedure used allows for the construction of 
lower bounds that are tighter and more broadly applicable than existing tail approximations. 
\end{abstract}

\maketitle

\section{Introduction}

This paper presents and solves a nonlinear optimization problem arising in the construction of lower bounds for the tails of distributions which possess a moment generating function on an open subset of $\R^+$.  The resulting lower bounds complement the classical Chernoff (upper) bound \cite{chernoff} in a set of cases more general than has been previously achieved.  

The methodology for the construction of the bounds was motivated by the presentation of the lower bound in Cram\'{e}r's large deviations theorem on p.29 of \cite{stroock}.  An earlier version of the results presented here was used in \cite{theodosopoulos} to establish a lower bound to the asymptotic convergence rate for an algorithm for global optimization.

The bounds presented here share numerous methodological characteristics with the development of {\it saddleppoint approximations} \cite{daniels1,daniels2,jensen}.  Both schemes use {\it tilting}, a technique first developed by Esscher, in order to center the power series expansions at the desired tail of the distribution.  Our method concentrates on the restriction of the Laplace transform on the real line.  A nonlinear optimization problem is constructed by adding two degrees of freedom to the tilting procedure.  This allows us to obtain tighter lower bounds, which hold even in cases where existing lower bounds break down.

The same direction was explored independently in \cite{bagdasarov} where a rough lower bound is computed using some rudiments of the methodology utilized here.  We parameterize the problem more efficiently, thereby arriving at a lower bound which possesses significantly better tightness.  The tools used in this paper are also very similar to those employed by Vinogradov (\cite{vinogradov}), in that they both explore beyond the realm of applicability of the Cram\'{e}r condition.  The main difference between our work and Vinogradov's lies in the different questions we ask.  Vinogradov assumes tail properties and extends classical large deviations results for sums of random variables under conditions not covered by classical techniques.  On the other hand, we are interested in inferring tail estimates under minimal conditions on the Laplace transform.  In that sense, despite the similarity in techniques with \cite{vinogradov}, our logic is more akin to that employed in \cite{bagdasarov} and \cite{jensen}.

Throughout the paper it is assumed that we have access to estimates of $\Xi ( \cdot)$, the Legendre dual of the {\it cumulant transform} (logarithm of the Laplace transform).  In the next section we introduce the new lower bound, represented as a two-dimensional constrained nonlinear optimization problem.  The third section offers comparisons with three alternative lower bounds.  Following that, we proceed to solve to solve the nonlinear optimization problem, thus arriving at efficient numerical estimates.  We include a figure which illustrates the comparison of the new lower bound with existing alternatives.

\section{A Lower Bound to Complement the Chernoff Bound}

Let $X$ be a real-valued, positive random variable on a probability space 
$({\mathcal X},\mu)$.  Assume that $X$ has exponential moments with respect 
to $\mu$, i.e. $g(\xi) \doteq \ee^\mu \left[e^{\xi X} \right] < \infty$
for $\xi$ in some open set $(-\infty, \xi^\ast)$, where $\xi^\ast > 0$.  
Let the {\it rate function} be given by the Legendre transform of the 
{\it cumulant} (i.e. the logarithm of the moment generating function),$I_\mu (y) \doteq \sup_\xi \{\xi y - \log g(\xi) \}$.
Further, let 
$$\Xi (y) \doteq \left\{ \begin{array}{ll} \arg \sup_{\xi \geq 0} \{\xi y - 
\log g(\xi) \} & \mbox{if $y \geq \ee^\mu [X]$} \\
\arg \sup_{\xi \leq 0} \{\xi y - \log g(\xi) \} & \mbox{if $y < \ee^\mu 
[X]$} \end{array} \right.$$
be the corresponding `Legendre dual'.  It is well known \cite{stroock} that 
$\Xi$ is an 
increasing concave function with $\Xi \left(\ee^\mu [X] \right) = 0$, and thus, 
it is generally invertible.  Using integration by parts we notice that
\begin{eqnarray}
& & {\frac {d}{dy}} g(\Xi(y)) = y g(\Xi(y)) \Xi'(y) \Longrightarrow 
\nonumber\\
& \Longrightarrow & g(\Xi(y)) = \exp \left\{y \Xi(y) - \int_{\ee^\mu [X]}^y
\Xi(t) dt \right\}. \label{eq:xirep}
\end{eqnarray}
The above formula leads to the following concise representation of the rate 
function:
\begin{equation}
I_\mu (y) =  \int_{\ee^\mu [X]}^y \Xi(t)dt.  \label{eq:rateint}
\end{equation}
Of course for any $\xi < \xi^\ast$,
$$\int_{\Xi^{-1} (\xi)}^x \Xi(t) dt = \int_\xi^{\Xi(x)} {\frac {\lambda 
d\lambda}{\Xi' (\Xi^{-1} (\lambda))}},$$
and
$$\Xi' \left(\Xi^{-1} (\lambda) \right) = {\frac {g^2}{g g'' 
- g'^2}} = \left( {\frac {d}{d\lambda}}{\frac {g'}{g}} 
\right)^{-1},$$
so that 
\begin{equation}
\int_{\Xi^{-1} (\xi)}^x \Xi(t) dt = x \Xi(x) -\xi \Xi^{-1} (\xi) -\log {\frac 
{g(\Xi(x))}{g(\xi)}}. \label{eq:nointegral}
\end{equation}
However, the integral representation of the exponent in (\ref{eq:xirep}) has 
two advantages.  Firstly it does not depend explicitly on the moment generating 
function, creating the possibility of constructing the tail bounds we are after 
using only the Legendre dual, $\Xi$, side-stepping the moment 
generating function.  We will explore this idea in a subsequent paper.  
Secondly, the integral representation of the exponent in (\ref{eq:xirep}) allows 
us to combine rate functions in a straightforward manner by adding the 
corresponding Legendre duals.

An application of the Markov inequality to the random variable $\exp\{\xi X\}$ 
suffices to obtain a surprisingly accurate upper bound to the tails of $X$, the 
celebrated Chernoff bound \cite{chernoff,bahadur} (assuming $y > 
\ee^\mu [X]$):

$$\mu(X \geq y) = \inf_{\xi \geq 0} \mu \left(e^{\xi X} \geq e^{\xi y} \right) 
\leq \inf_{\xi \geq 0} {\frac {g(\xi)}{e^{\xi y}}} \leq \exp \{-I_\mu (y)\}.$$

Note that this bound works equally well when $y < \ee^\mu [X]$, by considering the
left hand tail instead of the right hand one.  Specifically we obtain 
$$\mu (X \leq y) = \mu (-X \geq -y) = \inf_{\xi \leq 0} \mu \left(e^{\xi X} \geq 
e^{\xi y} \right),$$
which leads to the same upper bound for the left hand tail as the one for the 
right hand tail we obtained above in the case $y > \ee^\mu [X]$.

This short exposition of the Chernoff bound emphasizes the \textit{optimizing} 
degrees of freedom afforded to us by the free parameter $\xi$.  The estimates 
in the rest of the paper involve the construction of a lower bound to accompany 
the Chernoff bound. 

The tools we use in the construction of the desired lower bound are exponential 
tilting and the incorporation of a second optimizing degree of freedom.  The
former leads to a \textit{centering} of the measure around the tail of interest 
\cite{jensen}.  The latter allows us to tailor the tilting procedure in order 
to optimize the iterative use of the Chernoff bound to the tilted measure.

In particular, let $\nu_\xi \in \myprob ({\mathcal X})$ be a new probability 
measure on ${\mathcal X}$ defined by $\nu_\xi (dt) \doteq g(\xi)^{-1} e^{\xi t} \mu (dt)$.
This is called an exponentially tilted measure, after the theory of Esscher 
Tilting \cite{jensen}.  To simplify the notation, we will use $I_\alpha$ to 
signify $I_{\nu_{\Xi (\alpha y)}}$ and $\ee^\alpha$ for $\ee^{\nu_{\Xi (\alpha y)}}$.  
Observe that $I_\mu = I_{\scriptscriptstyle 
\ee^\mu [X]/y}$. 
Finally, let
$$L(\alpha, \delta, y) \doteq \left(1 - e^{-I_\alpha (\delta y)} - e^{- 
I_\alpha (y)} \right) \exp \left\{-I_\mu (\alpha y) - \Xi(\alpha y) y 
(\delta -\alpha) \right\}.$$
With this terminology we are in a position to state the proposed reversion of the Chernoff bound as a nonlinear constrained optimization in two dimensions:

\begin{theorem}
\label{theo1}
For any $y> \ee^\mu[X]$, the following inequality holds:
\begin{equation}
\mu(X \geq y) \geq \sup_{1<\alpha<\delta} L(\alpha, \delta, y). \label{eq:lb1}
\end{equation}
Moreover, there exist feasible values of $\alpha$ and $\delta$ which make the 
right hand side of (\ref{eq:lb1}) strictly positive.
\end{theorem}
\begin{proof}
Following the traditional proof of Cram\'{e}r's theorem \cite{stroock} we let $Y$ to 
be a $\nu_\alpha$-distributed random variable.  Observe that, for any $\alpha$,
$$\ee^\alpha [Y] = g \left( \Xi (\alpha y) \right)^{-1} \int_{-\infty}^\infty 
t e^{t \Xi (\alpha y)} \mu(dt) = \left. {\frac {d}{d\xi}} \right|_{\xi = \Xi(\alpha y)} 
\log g = \alpha y.$$
Then, for every $1<\alpha<\delta$ we have
\begin{eqnarray*}
\mu(X \geq y) & = & g \left( \Xi (\alpha y) \right) \int_y^\infty e^{-t 
\Xi (\alpha t)} \nu_\alpha (dt) \nonumber \\
& \geq & g \left( \Xi(\alpha y) \right) \int_y^{\delta y} 
e^{-t \Xi (\alpha t)} \nu_\alpha (dt) \nonumber \\
& \geq & g \left( \Xi(\alpha y) \right) e^{-\delta y \Xi (\alpha y)} 
\nu_\alpha \left( [y, \delta y] \right) \label{eq:lb1h}
\end{eqnarray*}
We are now in a position to apply the Chernoff bound iteratively as it were to 
estimate the last term on the right hand side.  Specifically, we observe that, for 
any $\delta > \alpha >1$, 
$\nu_\alpha (Y > \delta y) \leq \exp \left \{-I_\alpha (\delta y) \right\}$, because 
$\delta y > \ee^{\nu_\alpha} [Y] = \alpha y$  and $\nu_\alpha (Y < y) \leq \exp 
\left \{ -I_\alpha (y) \right\}$, because $y < \ee^{\nu_\alpha} [Y] = \alpha y$.  
Consequently we can estimate the last term on the right hand side of (\ref{eq:lb1h}) as
\begin{equation}
\nu_\alpha ([y, \delta y]) = 1 - \nu_\alpha (Y > \delta y) - \nu_\alpha (Y < y)
\geq 1 - e^{-I_\alpha (\delta y)} - e^{-I_\alpha (y)}  \label{eq:lb1hh}
\end{equation}
Substituting (\ref{eq:lb1hh}) into (\ref{eq:lb1h}) we see that, for any 
$1<\alpha <\delta$, $\mu (X \geq y) \geq L(\alpha, \delta)$.
Noting that the left hand side does not depend on $\alpha$ or $\delta$, we conclude 
that the inequality is maintained if we maximize the right hand side with respect to 
$\alpha$ and $\delta$, thus obtaining (\ref{eq:lb1}).

In order to evaluate the rate function for the tilted measure we observe that, for 
any $\theta$, 
$$\ee^\alpha \left[e^{\theta Y} \right] = {\frac {g \left( \theta + \Xi (\alpha y) 
\right)}{g \left( \Xi (\alpha y) \right)}}.$$
Thus, $\Xi_\alpha (\delta y) \doteq \arg \sup \left\{ \delta y \theta - \log 
\ee^\alpha \left[ e^{\theta Y} \right] \right\}$ must satisfy 
$$\delta y = {\frac {d}{d \theta}} \log \left( {\frac {g \left( \theta + \Xi (\alpha y) 
\right)}{g \left( \Xi (\alpha y) \right)}} \right) = \left. {\frac {d}{d\xi}} 
\right|_{\xi = \theta + \Xi (\alpha y)} \log g$$
and therefore
$$\Xi_\alpha (\delta y) = \Xi (\delta y) - \Xi (\alpha y).$$
Thus, using (\ref{eq:rateint}) in this situation we obtain
\begin{equation}
I_\alpha (\delta y) = \int_{\ee^\alpha [Y]}^{\delta y} \Xi_\alpha (t) dt = y (\alpha -\delta) \Xi (\alpha y) + \int_{\alpha y}^{\delta y} \Xi (t) dt 
\label{eq:ialphadeltay}
\end{equation}
and
\begin{equation}
I_\alpha (y) = \int_{\ee^\alpha [Y]}^y \Xi_\alpha (t) dt = y (\alpha -1) \Xi (\alpha y) - \int_y^{\alpha y} \Xi (t) dt \label{eq:ialphay}
\end{equation}
We are now in a position to show that there exists a feasible choice of $\alpha$ and 
$\delta$ such that $e^{-I_\alpha (\delta y)} + e^{-I_\alpha (y)} < 1$, thus ensuring 
that the right hand side of (\ref{eq:lb1}) is strictly positive.  Specifically, observe 
that the monotonicity of $\Xi (\cdot)$ leads to
\begin{eqnarray*}
{\frac {d}{d \delta}} I_\alpha (\delta y) & = & y \left( \Xi (\delta y) - \Xi (\alpha y)
\right) > 0 \\
{\frac {d^2}{d \delta^2}} I_\alpha (\delta y) & = & y^2 \Xi' (\delta y) \geq 0
\end{eqnarray*}
which implies that 
\begin{equation}
\lim_{\delta \rightarrow \infty} I_\alpha (\delta y) = \infty. \label{eq:limdelta}
\end{equation}
On the other hand, 
$${\frac {d}{d \alpha}} I_\alpha (y) = y^2 (\alpha -1) \Xi' (\alpha y) \geq 0$$
which, together with the observation that $I_{\alpha = 1} (y) = 0$, implies that there 
exists an $\epsilon>0$ and $1<\alpha^\ast <\infty$ such that $I_{\alpha^\ast} (y) 
=\epsilon$.  Choose a $\delta^\ast < \infty$ such that $I_\alpha (\delta y) > - \log
\left( 1- e^{-\epsilon} \right)$.  We can always do that because of 
(\ref{eq:limdelta}).  Then, manifestly, $I_{\alpha^\ast} (y) + I_{\alpha^\ast} 
(\delta y) < 1$ and therefore the right hand side of (\ref{eq:lb1}) is strictlly 
positive.
\end{proof}

Note that the lower bound in (\ref{eq:lb1}) does not depend on explicit knowledge of the 
moment generating function $g$.  Indeed, (\ref{eq:rateint}), (\ref{eq:ialphadeltay}) and 
(\ref{eq:ialphay}) show that all the components of (\ref{eq:lb1}) can be computed directly 
from $\Xi(\cdot)$.  This opens the possibility, discussed further in the Conclusions, that 
no precise knowledge, or perhaps even existence, of the moment generating function may be 
required for (\ref{eq:lb1}).

Furthermore, observe that the lower bound in (\ref{eq:lb1}) can be described without the 
use of any integrals.  Specifically, using (\ref{eq:nointegral}) with (\ref{eq:rateint}), 
(\ref{eq:ialphadeltay}) and (\ref{eq:ialphay}) we obtain 
\begin{eqnarray*}
I_\mu (\alpha y) & = & \alpha y \Xi(\alpha y) - \log g \left(\Xi (\alpha y) \right) \\
I_\alpha (\delta y) & = & \delta y \left (\Xi (\delta y) - \Xi (\alpha y) \right) + 
\log {\frac {g \left( \Xi (\alpha y) \right)}{g \left( \Xi ( \delta y) \right)}} \\ 
I_\alpha (y) & = & y \left (\Xi (y) - \Xi (\alpha y) \right) + 
\log {\frac {g \left( \Xi (\alpha y) \right)}{g \left( \Xi ( y) \right)}}
\end{eqnarray*}
Thus we see that all the components of $L(\alpha, \delta,y)$, and thus of (\ref{eq:lb1}), 
can be expressed without the need for any intergrals.  The integral representations shown 
above serve to do away with the explicit dependence on the moment generating function $g$.

\section{Comparison with Existing Lower Bounds}

At this point it is worthwhile to compare the lower bound in Theorem \ref{theo1} to 
three other approximations.  The first one is Daniels' saddlepoint approximation 
which, using our notation, is given by \cite{daniels1,daniels2}
\begin{equation}
\mu(X \geq y) \sim (2 \pi)^{-1/2} \int_y^\infty \sqrt{\Xi' (t)} e^{-I_\mu (t)} dt.
\label{eq:saddlepoint}
\end{equation}
Compared to (\ref{eq:lb1}) in Theorem \ref{theo1}, (\ref{eq:saddlepoint}) has the disadvantage that it involves an extra integral.  In that sense, the Chernoff bound and (\ref{eq:lb1}) can be thought as upper and lower bounds to the integral expression in (\ref{eq:saddlepoint}).

Second, we look at the lower bound proposed by Bagdasarov and Ostrovskii \cite{bagdasarov}. 
 While they use different notation, their methodology is very close to the one presented 
here.  Specifically, their lower bound to $\mu(X \geq y)$ has only one free parameter, 
$\Delta >0$, which, using the notation in the current paper, is equivalent to $1- 
\Xi(y)/\Xi(\alpha y)$.  As in the proof of Theorem \ref{theo1} above, their lower bound 
needs access to a point $y_+ > y$.  This point corresponds to the point $\delta y$ in our 
notation.  They describe this point as the point where the function $\lambda (1+\Delta) 
y I(y)$, where $\lambda$ takes the place of $\Xi (\alpha y)$ in the notation used 
here.  But the supremum of the function $\lambda (1+\Delta) y I(y)$ over $y$ is the 
Legendre transform of $I(y)$, which is itself the Legendre transform of $\log g(\xi)$.  
Using Legendre duality we conclude that $y_+ = \delta y = \Xi^{-1} \left( 2 \Xi(\alpha y) - 
\Xi (y) \right)$.  With these notational translations, the Bagdasarov-Ostrovskii (B-O) 
lower bound can be described as
\begin{equation}
\mu(X \geq y) \geq \sup_{\alpha > 1} \tilde{L} (\alpha, y), \label{eq:BOlb}
\end{equation}
where
$$\tilde{L} (\alpha, y) \doteq {\frac {1 - {\frac {\Xi(\alpha y)}{\Xi(\alpha y) - \Xi(y)}} 
\left[ e^{-I_\alpha \left( \delta(\alpha, y) y \right)} + e^{- 
I_\alpha (y)} \right]}{ 1 - e^{-I_\alpha \left( \delta(\alpha, y) y \right)} - e^{- I_\alpha
(y)}}} L \left( \alpha, \delta (\alpha, y), y \right),$$
and
$$\delta (\alpha, y) \doteq y^{-1} \Xi^{-1} \left( 2\Xi(\alpha y) - \Xi(y) \right).$$
So we can see immediately that the B-O lower bound is inferior to the one described in 
Theorem \ref{theo1} for two reasons.  On the one hand it foregoes one of the two optimizing 
degrees of freedom (making $\delta$ a function of $\alpha$, which we will recognize as a suboptimal choice in the following section).  Furthermore, the term ${\frac 
{\Xi(\alpha y)}{\Xi(\alpha y) - \Xi(y)}}$ makes the fraction on the right hand side of the 
expression for $\tilde{L}$ strictly less than 1.  

Also, \cite{bagdasarov} does not provide a general statement about the range of 
applicability of the B-O lower bound as Theorem \ref{theo1} does.  It turns out that there are cases of interest where the B-O lower bound is inapplicable.   While the B-O lower bound is less tight than (\ref{eq:lb1}) and its range of applicability is not as broad, \cite{bagdasarov} has the advantage of presenting their 
lower bound assuming only approximate knowledge of the moment generating function $g(\xi)$. 
 By contrast, the current paper assumes that we have complete knowledge of the moment 
generating function.  It turns out that this is not necessary.  Motivated by Bagdasarov and 
Ostrovskii's work we extend the results presented here to the more general case of only 
approximate knowledge of the moment generating function in a follow-up paper.

In the same spirit is the lower bound presented in \cite{stroock}.  The construction of 
Stroock's lower bound is very similar to the one we present here, and in fact our 
presentation mirrors his.  The main difference lies with Stroock's use of the 
Chebyshev inequality to bound $\nu_\alpha ([y, \delta y])$, as opposed to our iterative use 
of the Chernoff bound.  The symmetry of the Chebyshev inequality around the mean determines 
one of the two optimizing degrees of freedom, and consequently Stroock's lower bound can 
be described as:
\begin{equation}
\mu(X \geq y) \geq \sup_{\alpha > 1} \hat{L} (\alpha, y), \label{eq:stroock}
\end{equation}
where
$$\hat{L} (\alpha, y) \doteq {\frac {1 - {\frac {1}{\Xi^\prime (\alpha y) y^2 (\alpha 
-1)^2}}}{ 1 - e^{-I_\alpha \left((2 \alpha -1) y \right)} - e^{- I_\alpha (y)}}} L \left( 
\alpha, 2\alpha -1, y \right).$$
The first disadvantage of (\ref{eq:stroock}) when compared to (\ref{eq:lb1}) is the fact 
that it lacks one degree of freedom, whose optimization could only improve the latter.  
Secondly, unlike (\ref{eq:lb1} which is guarranteed to work in general by Theorem 
\ref{theo1}, the range of applicability of (\ref{eq:stroock}) is limited by the requirement 
that $\Xi^\prime (\alpha y) y^2 (\alpha -1)^2 >1$.  There are indeed application of 
interest (which will be discussed in a subsequent section) that do not conform with this 
requirement, and for which therefore (\ref{eq:stroock}) is inapplicable.  In particular, 
one such application involves $\Xi(t) = c - t^{-1}$ for some constant $c>0$.  One readily 
concludes that this choice makes $\Xi^\prime (\alpha y) y^2 (\alpha -1)^2 = \left(1 - 
{\frac {1}{\alpha}} \right)^2 <1$, thus invalidating (\ref{eq:stroock}).

Finally, even in its range of applicability, (\ref{eq:stroock}) is less tight than 
(\ref{eq:lb1}).  In order to see this, let's consider the first order approximation to 
$\Xi(\cdot)$ around $\alpha y$ we have
\begin{equation}
\Xi(t) = \Xi(\alpha y) + \Xi^\prime (\alpha y) (t- \alpha y), \label{eq:linearapprox}
\end{equation}
and thus, using (\ref{eq:ialphadeltay}, \ref{eq:ialphay}) we obtain
\begin{equation}
I_\alpha (y) = - \Xi^\prime (\alpha y) \left\{ {\frac {\alpha^2 y^2 - y^2}{2}} - 
\alpha y^2 (\alpha -1) \right\} = {\frac {\Xi^\prime (\alpha y) y^2 (\alpha -1)^2}{2}} \label{eq:linearialphay}
\end{equation}
and
\begin{eqnarray}
I_\alpha \left( (2\alpha -1) y \right) & = & \Xi^\prime (\alpha y) \left\{ {\frac {(2\alpha 
-1)^2 y^2 - \alpha^2 y^2}{2}} - \alpha y^2 (2\alpha -1 -\alpha) \right\} \nonumber \\
& = & {\frac {\Xi^\prime (\alpha y) y^2 (\alpha -1)^2}{2}} \label{eq:linearialphadeltay}
\end{eqnarray}
The following lemma shows that, under the linear approximation to $\Xi(\cdot)$, when 
(\ref{eq:stroock}) is valid, the ratio on the right hand side of the expression for 
$\hat{L}$ is less than 1 for $y$ large enough.
\begin{lemma}
\label{lemma1}
Fix $\alpha >1$ and $y$ such that $I_\alpha (y) > {\frac {1}{2}}$.  Assume that 
$\Xi(\cdot)$ is a linear function.  Then, for $y$ large enough,
\begin{equation}
\Xi^\prime (\alpha y) y^2 (\alpha -1)^2 \left[ e^{-I_\alpha \left( (2\alpha -1) y 
\right)} + e^{-I_\alpha (y)} \right] < 1. \label{eq:lemma1}
\end{equation}
\end{lemma}
\begin{proof}
Using (\ref{eq:linearialphay}) and (\ref{eq:linearialphadeltay}) we can rewrite the 
expression on the left hand side of (\ref{eq:lemma1}) as $4I_\alpha (y) e^{-I_\alpha (y)}$. 
Consider the function $4we^{-w}$; it is clear that for $w$ large enough (in particular 
$w>2.16$ would suffice), $4we^{-w}<1$.

Observe that ${\frac {d}{dt}} I_\alpha (t)= \Xi(t) - \Xi(\alpha y)$, which is zero only at 
$t=\alpha y$.  This unique critical point is a minimum since $\left. {\frac {d^2}{dt^2}} \right|_{t=\alpha y} I_\alpha = \Xi^\prime (\alpha y) 
\geq 0$,
because $\Xi(\cdot)$ is a concave non-decreasing function \cite{stroock}.  Therefore, 
$\inf_t I_\alpha (t) = I_\alpha (\alpha y) = \alpha y (\alpha -1) \Xi (\alpha y)$, which 
is monotonically increasing with $y$.  Thus $\lim_{y \rightarrow \infty} \inf_t I_\alpha (t) \leq \lim_{y \rightarrow \infty} 
I_\alpha (y) = \infty$.

Putting the last two statements together we see that, indeed, for large enough $y$, the 
left hand side of (\ref{eq:lemma1}) is strictly less than 1.
\end{proof}
Lemma \ref{lemma1} immediately implies that, even when (\ref{eq:stroock}) is applicable, 
it is less tight than (\ref{eq:lb1}) far enough in the tail.  

\section{A Nonlinear Optimization Problem}

We are now in a position to compute the lower bound presented in an implicit way in 
(\ref{eq:lb1}).  In order to arrive at an explicit computation, we need to solve the 
optimization over the two parameters, $\alpha$ and $\delta$, that determine the tightest 
achievable lower bound.

The first step in our computation is the reduction of the optimization in (\ref{eq:lb1}) to one variable.  In what follows we will use the following symbols to simplify the presentation:

$$A(\alpha, y) \doteq \exp\left\{ - I_\alpha (y) \right\}$$
$$B(\alpha, \delta, y) \doteq \exp\left\{ - I_\alpha (\delta y) \right\}$$

We proceed by evaluating the first order condition with respect to $\alpha$:

\begin{eqnarray}
{\frac {\partial L}{\partial \alpha}} (\alpha, \delta, y) & = & -\delta y^2 \Xi' (\alpha y) L + \alpha y^2 \Xi' (\alpha y) L -{\frac {Ly^2 \Xi' (\alpha y)}{1-A-B}} \left\{ (\delta -\alpha)B -(\alpha -1) A \right\} = 0 \nonumber\\
\Longleftrightarrow & & -\delta + \alpha -{\frac {(1-\alpha)A + (\delta -\alpha)B}{1-A-B}} = 0 \Longleftrightarrow \nonumber\\
& & \delta^\ast (\alpha, y) = {\frac {\alpha -A}{1-A}} \label{eq:deltastar}
\end{eqnarray} 
The following lemma describes the properties of the resulting optimum choice of $\delta$:
\begin{lemma}
\label{lemma2}
For every $y$, $\delta^\ast$ is a quasiconvex function of $\alpha$ which attains a unique minimum at some $\check{\alpha} (y) \in (1, \infty)$.
\end{lemma}
\begin{proof}
Observe that 
$$\lim_{\alpha \rightarrow 1^+} A = \exp\left\{ \lim_{\alpha \rightarrow 1^+} \left[ -y(\alpha -1) \Xi(\alpha y) + \int^{\alpha y}_y \Xi(t) dt \right] \right\} =  1^- $$
and
\begin{equation}
A_\alpha \doteq {\frac {\partial A}{\partial \alpha}} = -Ay^2(\alpha -1) \Xi' (\alpha y). \label{eq:dAdalpha}
\end{equation} 
Also, $\lim_{\alpha \rightarrow \infty} A = 0^+$ because the concavity of $\Xi$ \cite{stroock} implies $\lim_{\alpha \rightarrow \infty} \Xi' (\alpha y) < \infty$ and therefore $\lim_{\alpha \rightarrow \infty} {\frac {\partial}{\partial \alpha}} \log A = -\infty$, which implies that $\log A$, and therefore $A$ approach 0 from above as $\alpha$ tends to infinity.  
The proof of the lemma will take three steps.  The first step is to show that, for any $y$, 
\begin{equation}
\lim_{\alpha \rightarrow 1^+} \delta^\ast = \lim_{\alpha \rightarrow 1^+} {\frac {1 -A_\alpha}{-A_\alpha}} = \lim_{\alpha \rightarrow 1^+} {\frac {1}{Ay^2 (\alpha -1) \Xi' (\alpha y)}} = + \infty. \label{eq:deltastarof1}
\end{equation}
The second step is to observe that $\lim_{\alpha \rightarrow \infty} \delta^\ast = +\infty$.
The final step of the proof involves
\begin{equation}
{\frac {\partial \delta^\ast}{\partial \alpha}} = {\frac {1 -A -Ay^2 (\alpha -1)^2 \Xi' (\alpha y)}{(1-A)^2}}. \label{eq:ddeltastardalpha}
\end{equation}
From (\ref{eq:deltastarof1}) we conclude that $\lim_{\alpha \rightarrow 1^+} {\frac {\partial \delta^\ast}{\partial \alpha}} = - \infty$.  Therefore, for every $y$, there must exist a $\check{\alpha} (y) \in (1, \infty)$ such that ${\frac {\partial \delta^\ast}{\partial \alpha}} \left(\check{\alpha} (y), y \right) = 0$
and 
$\forall \alpha \in \left( \check{\alpha} (y), \infty \right), {\frac {\partial \delta^\ast}{\partial \alpha}} \left(\check{\alpha} (y), y \right) > 0$.
Indeed, differentiating (\ref{eq:ddeltastardalpha}) with respect to $\alpha$ we see that
\begin{eqnarray}
{\frac {\partial^2 \delta^\ast}{\partial \alpha^2}} & = & {\frac {Ay^2 (\alpha -1)}{(1-A)^2}} \left[y^2 (\alpha-1)^2 \Xi' (\alpha y)^2 + \Xi' (\alpha y) -y (\alpha-1) \Xi'' (\alpha y) \right] + \nonumber \\
& & + {\frac {2Ay^2 (\alpha-1) \Xi' (\alpha y)}{1-A}} {\frac {\partial \delta^\ast}{\partial \alpha}}. \label{eq:dsqdeltastardalphasq}
\end{eqnarray}
The concavity of $\Xi$ together with (\ref{eq:dsqdeltastardalphasq}) show us that ${\frac {\partial \delta^\ast}{\partial \alpha}} = 0 \Longrightarrow {\frac {\partial^2 \delta^\ast}{\partial \alpha^2}} > 0$.
Implying that each critical point of $\delta^\ast$ is a minimum, and therefore there is a unique minimum and $\delta^\ast$ is quasiconvex.
\end{proof}

\begin{figure}
\epsfxsize=3.8in
\epsfbox{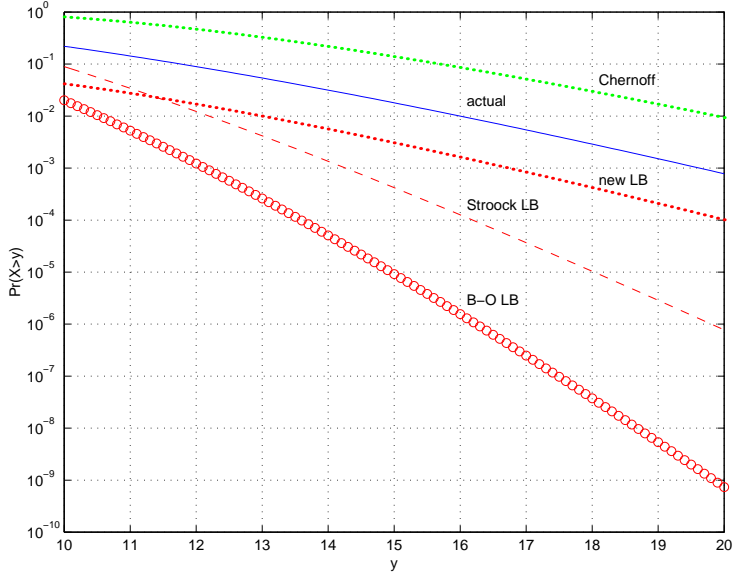}
\caption{The right tail of a gamma variate with pdf given by $f(x)={\frac {x^7 e^{-x}}{\Gamma (8)}}$.}
\label{fig:gammatails}
\end{figure}

Using the resulting expression for $\delta^\ast$ as a function of $\alpha$, we can rewrite the lower bound $L$ and the expression $B$ above as functions solely of $\alpha$ and $y$:
$$\hat{L} (\alpha, y) \doteq L \left( \alpha, \delta^\ast (\alpha, y), y \right)$$
$$\hat{B} (\alpha, y) \doteq \exp\left\{ - I_\alpha \left(y\delta^\ast (\alpha, y) \right) \right\}$$
Using this terminology, we observe that:
\begin{lemma}
\label{lemma3}
For every $y$ there exists a unique $\hat{\alpha} (y) \in (1, \check{\alpha} (y)$ such that $A(\hat{\alpha} (y), y) + \hat{B} (\hat{\alpha}, y) = 1$ and, for all $\alpha \in (1, \hat{\alpha})$, $A(\alpha, y) + \hat{B} (\alpha, y) < 1$.
\end{lemma}
\begin{proof}
Notice that
\begin{equation}
\hat{B}_\alpha \doteq {\frac {\partial \hat{B}}{\partial \alpha}} = \hat{B}y \left\{y \left(\delta^\ast -\alpha \right) \Xi^\prime (\alpha y) -{\frac {\partial \delta^\ast}{\partial \alpha}} \left[ \Xi (\delta^\ast y) -\Xi (\alpha y) \right] \right\} \label{eq:dBhatdalpha}
\end{equation}
By the concavity of $\Xi$ and Lemma \ref{lemma2} we can see that $\hat{B}_\alpha \geq 0$ and therefore 
\begin{equation}
\lim_{\alpha \rightarrow \infty} \hat{B} >0. \label{eq:Bhatofinfty}
\end{equation}
Using (\ref{eq:dAdalpha}) and (\ref{eq:dBhatdalpha}), after some algebra we arrive at
\begin{equation}
{\frac {\partial}{\partial \alpha}} (A + \hat{B}) = y^2 \Xi^\prime (\alpha y) (\delta^\ast -\alpha) \left\{A+B-1 -{\frac {\partial \delta^\ast}{\partial \alpha}} {\frac {\Xi(\delta^\ast y) -\Xi(\alpha y)}{y(\delta^\ast -\alpha) \Xi^\prime (\alpha y)}} \right\}. \label{eq:dAplusBhatdalpha}
\end{equation}
Since the term outside the brackets on the right hand side of (\ref{eq:dAplusBhatdalpha}) is always nonnegative, we conclude, using Lemma \ref{lemma2}, that the following statements are true:
\begin{enumerate}
\item $\exists \bar{\alpha} \in (1, \check{\alpha}) \left. {\frac {\partial}{\partial \alpha}} \right|_{\bar{\alpha}} (A + \hat{B}) = 0 \Longrightarrow A(\bar{\alpha}) + \hat{B} (\bar{\alpha}) <1$.
\item $\exists \bar{\alpha} \in (\check{\alpha}, \infty) \left. {\frac {\partial}{\partial \alpha}} \right|_{\bar{\alpha}} (A + \hat{B}) = 0 \Longrightarrow A(\bar{\alpha}) + \hat{B} (\bar{\alpha}) >1$.
\item $\exists \alpha^\ast \in (1, \check{\alpha}) A(\alpha^\ast) + \hat{B} (\alpha^\ast) = 1 \Longrightarrow \left. {\frac {\partial}{\partial \alpha}} \right|_{\alpha^\ast} (A + \hat{B}) > 0$.
\item $\exists \alpha^\ast \in (\check{\alpha}, \infty) A(\alpha^\ast) + \hat{B} (\alpha^\ast) = 1 \Longrightarrow \left. {\frac {\partial}{\partial \alpha}} \right|_{\alpha^\ast} (A + \hat{B}) < 0$.
\end{enumerate}
Also, using (\ref{eq:deltastarof1}), we see that
\begin{eqnarray*}
\lim_{\alpha \rightarrow 1^+} \hat{B} & = & \exp \left\{ \lim_{\delta \rightarrow \infty} \delta \left[ y \Xi(y) \left(1 -{\frac {1}{\delta}} \right) -{\frac {\int_y^{\delta y} \Xi(t) dt}{\delta}} \right] \right\} \\
& = & \exp \left\{ \lim_{\delta \rightarrow \infty} \delta y \left[ \Xi(y) -\Xi (\delta y) \right] \right\}= 0^+
\end{eqnarray*}
whether $\lim_{y \rightarrow \infty} \Xi < \infty$ or $\lim_{y \rightarrow \infty} \Xi = \infty$.  This deduction rules out a maximum of $A+\hat{B}$ before $\check{\alpha}$, because by statement (1) above, any critical point of $A+\hat{B}$ before $\check{\alpha}$ leads to $A+\hat{B}<1$ and therefore must be a minimum. 
Furthermore, the combination of statements (2) and (4) above imply that, if $A+\hat{B}$ lacks a zero below $\check{\alpha}$, then it cannot have a minimum below $\check{\alpha}$, because when $A +\hat{B}<1$ for all $\alpha < \check{\alpha}$, the slope of $A +\hat{B}$ will be negative for all $\alpha < \check{\alpha}$.  Therefore, one of the following two statements must hold:
\begin{enumerate}
\item[(i)] Either $A+\hat{B}>1$ for all $y$ and $\alpha >1$, or
\item[(ii)] There exists a zero of $A+\hat{B}$, $\hat{\alpha}$, in $(1, \check{\alpha})$ such that $\left.{\frac {\partial}{\partial \alpha}}\right|_{\hat{\alpha}} (A+\hat{B})>0$.   
\end{enumerate}
It turns out that we can rule out case (i).  Specifically, for every $\alpha$ and $y$, let $\check{\delta} (\alpha, y) \geq \alpha$ be such that $A(\alpha, y) + B\left( \alpha, \check{\delta} (\alpha, y), y \right) = 1$.  Clearly, 
$${\frac {\partial B}{\partial \delta}} = -By \left\{ \Xi (\delta y) - \Xi (\alpha y) \right\} \leq 0.$$
Thus, if there is no $\alpha \in (1, \check{\alpha})$ with $A+\hat{B}<1$, then for all $\alpha \in (1, \check{\alpha})$ and every $y$, $\delta^\ast (\alpha, y) < \check{\delta} (\alpha, y)$, which implies that
$$L_\alpha \doteq {\frac {\partial L}{\partial \alpha}} = {\frac {L y^2 \Xi^\prime (\alpha y)}{1-A-B}} \left[ (\alpha -A) -\delta (1 -A) \right] \leq 0.$$
But this would imply that the maximum lower bound $L$ is achieved on $\check{\delta}$, which leads to $\max L = 0$.  This clearly contradicts the statement of Theorem \ref{theo1} which asserts that, for any $y$, there exists a pair of values for $\alpha$ and $\delta$ making $L$ strictly positive.  Thus, we are left with statement (ii) as the only viable possibility, which establishes the desired result. 
\end{proof}
Let $G(\alpha, \delta, y) \doteq B(\alpha, \delta, y) \Xi (\delta y) -\left( 1 -A(\alpha, y) \right) \Xi(\alpha y)$.  
\begin{lemma}
\label{lemma4}
For any $\alpha>1$ and $y$, there exists a unique $\hat{\delta}(\alpha, y) \in (\alpha, \infty)$ such that $G\left(\alpha, \hat{\delta} (\alpha, y), y \right)=0$.  Moreover, ${\frac {\partial G}{\partial \delta}} \left(\alpha, \hat{\delta} (\alpha, y), y \right) <0$. 
\end{lemma}
\begin{proof}
From the definition of $G$ we see that
$$G_\delta \doteq {\frac {\partial G}{\partial \delta}} = -By \left\{\Xi (\delta y) \left[ \Xi(\delta y) -\Xi(\alpha y) \right] -\Xi^\prime (\delta y) \right\}$$
and 
$$G_{\delta \delta} \doteq {\frac {\partial^2 G}{\partial \delta^2}} = -y \left[ \Xi(\delta y) -\Xi(\alpha y) \right] {\frac {\partial G}{\partial \delta}} -By^2 \left\{ \Xi^\prime (\delta y) \left[ 2\Xi (\delta y) - \Xi (\alpha y) \right] - \Xi'' (\delta y) \right\},$$
which implies that $G_\delta \geq 0 \Longrightarrow G_{\delta \delta} <0$.  Observe that $G_\delta (\alpha, \alpha,y) = By \Xi^\prime (\alpha y) >0$.  Also, using (\ref{eq:Bhatofinfty}) we can see that for large enough $\delta$, $G_\delta (\alpha, \delta, y)<0$, because the concavity of $\Xi$ forces $\Xi^\prime$ to remain uniformly bounded.  Thus, for any $\alpha>1$ and $y$, $G(\alpha,\cdot, y)$ has a unique maximum $\bar{\delta} (\alpha, y) \in (\alpha, \infty)$.  
Notice that $G(\alpha, \alpha, y) = A \Xi(\alpha y)$ because $B(\alpha, \alpha, y) = 1$.  Also, for any $\alpha >1$ and $y$, 
\begin{eqnarray*}
\lim_{\delta \rightarrow \infty} B(\alpha, \delta, y) & = & \exp \left\{ \lim_{\delta \rightarrow \infty} \delta \left[ y \left( 1 –{\frac {\alpha}{\delta}} \right) \Xi(\alpha y) –{\frac {\int_{\alpha y}^{\delta y} \Xi(t) dt}{\delta}} \right] \right\} \\
& = & \exp \left\{ \lim_{\delta \rightarrow \infty} \delta y \left[ \Xi(\alpha y) -\Xi(\delta y) \right] \right\}= 0^+
\end{eqnarray*}
whether $\lim_{y \rightarrow \infty} \Xi < \infty$ or $\lim_{y \rightarrow \infty} \Xi = \infty$.  Therefore, $G$ indeed possesses a zero, $\hat{\delta} (\alpha ,y)$.  Naturally, $\hat{\delta} > \bar{\delta}$, and therefore ${\frac {\partial G}{\partial \delta}} \left(\alpha, \hat{\delta}(\alpha, y), y \right) <0$.  Finally, this zero is unique because for there to be another, there first must exist a minimum, which is prohibited by the preceding.
\end{proof}
We are in a position to prove the main theorem of this section:

\begin{theorem}
\label{theo2}
For every $y$ there exists a $\alpha^\ast (y) \in \left(1, \hat{\alpha} (y) \right)$ which attains the unique maximum of $\hat{L}$.
\end{theorem}
\begin{proof}
We have already seen that, for any $\alpha >1$ and any $y$, $L_\alpha \left(\alpha, \delta^\ast (\alpha, y), y \right) = 0$.  We can also see that
$$L_\delta \doteq {\frac {\partial L}{\partial \delta}} = {\frac {LGy}{1-A-B}}$$
and thus for any $\alpha >1$ and any $y$, $L_\delta \left(\alpha, \hat{\delta} (\alpha, y), y \right) = 0$.  The discussion at the end of the proof of Lemma \ref{lemma3} guarantees, for any $y$, the existence of an intersection, $\alpha^\ast \in (1, \hat{\alpha} (y))$ between the two curves, $\delta^\ast (\cdot, y)$ and $\hat{\delta} (\cdot, y)$.  The only question that remains is the uniqueness of this intersection and consequently of the maximum for $L$.  
Observe that, for any $y$,
\begin{eqnarray}
& & G \left(\alpha, \hat{\delta} \right) = 0 \Longrightarrow B \left( \alpha, \hat{\delta} \right) \Xi \left( y \hat{\delta} \right) = (1 -A) \Xi (\alpha y)  \nonumber \\ 
& \Longrightarrow & \left[ {\frac {\partial B}{\partial \alpha}} \left( \alpha, \hat{\delta} \right) + {\frac {\partial B}{\partial \delta}} \left(\alpha, \hat{\delta} \right) {\frac {\partial \delta^\ast}{\partial \alpha}} \right] \Xi \left( y \hat{\delta} \right) + B \left(\alpha, \hat{\delta} \right) y {\frac {\partial \delta^\ast}{\partial \alpha}}  \Xi^\prime \left( y \delta^\ast  \right) \nonumber \\
& & = -{\frac {\partial A}{\partial \alpha}} \Xi (\alpha y) = y (1 -A) \Xi^\prime (\alpha y) \nonumber \\
& \Longrightarrow & {\frac {\partial \hat{\delta}}{\partial \alpha}} = {\frac {y^2 \Xi (\alpha y) \Xi^\prime (\alpha y) \left[ A(\alpha -1) –(1-A)(\hat{delta} -\alpha) \right] + y \Xi^\prime (\alpha y) (1-A)}{B \left(\alpha, \hat{\delta} \right) \left[ \Xi \left(y \hat{\delta} \right) + y \Xi^\prime \left(y \hat{\delta} \right) \right]}} \label{eq:ddeltahatdalpha}
\end{eqnarray}
where the second line arises from differentiating both sides of the first line.  At the intersection $\alpha^\ast$ of $\delta^\ast$ and $\hat{\delta}$ we can simplify (\ref{eq:ddeltahatdalpha}) and obtain
$${\frac {y {\frac {\Xi^\prime (\alpha y)}{\Xi (\alpha y)}}}{1+ y {\frac {\Xi^\prime \left( y \hat{\delta} (\alpha, y) \right)}{\Xi \left( y \hat{\delta} (\alpha, y) \right)}}}} >0.$$
On the other hand, since $\alpha^\ast < \hat{\alpha} < \check{\alpha}$, we know that ${\frac {\partial \delta^\ast}{\partial \alpha}} < 0$.  Therefore, 
$$\delta^\ast = \hat{\delta} \Longrightarrow {\frac {\partial \delta^\ast}{\partial \alpha}} < {\frac {\partial \hat{\delta}}{\partial \alpha}}$$
and therefore the intersection of $\delta^\ast$ and $\hat{\delta}$ must be unique.
\end{proof}
Figure \ref{fig:gammatails} shows the lower bound resulting from the procedure described in this paper, compared to the alternatives discussed above, the Chernoff (upper) bound and the exact tail, which can be readily computed in this case.  This example was chosen in the range where the alternative lower bounds are also applicable, to allow for a comparison.  As we saw above, the new lower bound does not have the limitations in its range of applicability that plague the alternative lower bounds.  We can in see in Figure \ref{fig:gammatails} that the new lower bound maintains a consistent gap from Chernoff bound and exact tail, unlike the alternative lower bounds.
\vspace{0.2in}
\section{Conclusions and Future Steps}

We have shown a way to construct a lower bound to complement the Chernoff bound that is applicable generically, without restrictions to the moment generating function that characterized earlier inequalities of a similar nature.  We were able to represent the bound as the solution of a two-dimensional nonlinear optimization problem, and proved the existence and uniqueness of the solution.

As we saw earlier, the new lower bound has two technical advantages, aside from its broad applicability.  Specifically, unlike the saddlepoint approximation, it can be formulated without the need for any integrals.  Alternatively, it can be formulated using only $\Xi(\cdot)$, which makes it easy to combine across iid sequences.  From a theoretical perspective, the main advantage of the new lower bound is that it does not depend on an appeal to the law of large numbers, which often does not hold in situations of interest 

At this point it is natural to investigate the asymptotic properties (in the spirit of \cite{borovkov}) of the new lower bound, including its asymptotic gap from the Chernoff bound.  It is reasonable to expect that one can classify moment generating functions relative to the resulting asymptotic gaps.  While there are obvious examples with no gap (e.g. Gaussian) and some with a gap \cite{theodosopoulos}, a complete classification is still out of reach.  Extensions to the multivariate case are another natural next step, along the lines of \cite{barndorff,iyengar} for the saddlepoint approximation.

Furthermore, one may inquire whether the new lower bound, as well as the Chernoff bound, can be extended in situations where the moment generating function is not precisely known \cite{bagdasarov} or does not exist at all.  The former case will be dealt with in a follow-up paper.  More generally, the fact that, as shown in this paper, bilateral tail bounds are achievable with reference only to $\Xi(\cdot)$ lends support to the latter possibility.  It is also plausible to substitute the exponential function for other ones, more appropriate for different distributions.  An investigation of this question remains open at the moment.

\bibliographystyle{amsalpha}

\end{document}